\documentclass[pdflatex,sn-mathphys-num]{sn-jnl}

\usepackage{graphicx}
\usepackage{multirow}
\usepackage{amsthm}
\usepackage{mathrsfs}
\usepackage[title]{appendix}
\usepackage{textcomp}
\usepackage{manyfoot}
\usepackage{booktabs}
\usepackage{algorithmicx}
\usepackage{listings}
\usepackage{graphicx} 
\usepackage{amsmath}
\usepackage{amsfonts}
\usepackage{amssymb}
\usepackage{comment}
\usepackage{algorithm}
\usepackage{algpseudocode}
\usepackage[dvipsnames]{xcolor}
\usepackage{booktabs}
\usepackage{adjustbox}  
\usepackage{soul}

\theoremstyle{thmstyleone}

\theoremstyle{thmstyletwo}

\theoremstyle{thmstylethree}

\begin{document}

\title[Article Title]{\centerline{Neural Preconditioning via Krylov Subspace Geometry}}

\author[1]{\fnm{Nunzio} \sur{Dimola}}\email{nunzio.dimola@polimi.it}

\author*[2]{\fnm{Alessandro} \sur{Coclite}}\email{alessandro.coclite@poliba.it}

\author[1]{\fnm{Paolo} \sur{Zunino}}\email{paolo.zunino@polimi.it}

\affil[1]{\orgdiv{MOX, Dipartimento di Matematica}, \orgname{Politecnico di Milano}, \orgaddress{\street{Piazza Leonardo da Vinci 32}, \city{Milano}, \postcode{20133}, \country{Italy}}}

\affil[2]{\orgdiv{Dipartimento di Ingegneria Elettrica e dell'Informazione (DEI)}, \orgname{Politecnico di Bari}, \orgaddress{\street{Via Re David 200}, \city{Bari}, \postcode{70125}, \country{Italy}}}

\maketitle

\begin{abstract}

We propose a geometry-aware strategy for training neural preconditioners tailored to parametrized linear systems arising from the discretization of mixed-dimensional partial differential equations (PDEs). These systems are typically ill-conditioned because of the presence of embedded lower-dimensional structures and are solved using Krylov subspace methods. Our approach yields an approximation of the inverse operator employing a learning algorithm consisting of a two-stage training framework: an initial static pretraining phase, based on residual minimization, followed by a dynamic fine-tuning phase that incorporates solver convergence dynamics into training via a novel loss functional. This dynamic loss is defined by the principal angles between the residuals and the Krylov subspaces. It is evaluated using a differentiable implementation of the Flexible GMRES algorithm, which enables backpropagation through both the Arnoldi process and Givens rotations. The resulting neural preconditioner is explicitly optimized to improve early-stage convergence and reduce iteration counts in a family of 3D-1D mixed-dimensional problems with geometric variability of the 1D domain. Numerical experiments show that our solver-aligned approach significantly improves convergence rate, robustness, and generalization.


\end{abstract}
  
\section{Introduction}

The efficient numerical solution of parameter-dependent linear systems arising from discretized partial differential equations (PDEs) is a critical computational challenge in scientific computing. This task becomes particularly demanding in multi-query applications, such as uncertainty quantification or control problems, where the underlying PDE must be solved repeatedly for different parameter values. In this context, mixed-dimensional PDEs~\cite{D'Angelo20081481,10.1007/s10231-020-01013-1, Kuchta2021558, Heltai20232425} (just to make a few examples), which describe coupled operators defined on domains of varying dimensions, present additional challenges due to their inherent ill-conditioning, complex geometries, and heterogeneous material properties.
To address these large, sparse linear systems, iterative solvers, notably Krylov subspace methods such as the Conjugate Gradient (CG) and Generalized Minimal Residual (GMRES) algorithm, are widely employed~\cite{saad2003iterative}. However, the efficiency of these methods strongly depends on the spectral properties and condition number of the system matrices, which often deteriorate with increasing problem size or parameter variations. This fundamental limitation necessitates the extensive use of preconditioning techniques~\cite{https://doi.org/10.1002/nla.716}, which transform the original system into an equivalent one more amenable to rapid convergence. 

In this work, we focus on a family of linear systems arising from the discretizations of parameter-dependent mixed-dimensional PDEs
\begin{equation*}
\mathsf{ A^\mu u^\mu = b^\mu}, \quad \mu \in \mathcal{P},
\end{equation*}
where $\mathcal{P}$ denotes the parameter space encoding geometric variations of embedded lower-dimensional structures within a three-dimensional domain. The central challenge addressed here is the design of effective, parameter-dependent preconditioners that can handle significant geometric variability while ensuring robust solver convergence.
Traditional preconditioning strategies, including stationary methods (e.g., Jacobi, Gauss-Seidel), incomplete factorizations, sparse approximate inverses (SPAI) \cite{chen2005matrix}, multigrid (MG) \cite{trottenberg2000multigrid}, and domain decomposition (DD) methods \cite{quarteroni1999domain}, have been developed to improve convergence. For mixed-dimensional problems, physics-based block and algebraic multigrid methods tailored for metric-perturbed coupled problems have proven instrumental for efficient resolution~\cite{Kuchta2016B962, budisa2024algebraic, FIRMBACH2024117256}.
More recently, neural network-based strategies have gained significant attention for their ability to learn effective preconditioners directly from data~\cite{Azulay2023S127,doi:10.1137/24M162861X}. 
These approaches often leverage operator learning frameworks, like for example DeepONet~\cite{LuJinPang}, to approximate inverse operators or to construct transfer operators, or utilize the spectral bias of neural networks \cite{xu2024overview} to efficiently address low-frequency error components.

Our prior work~\cite{dimola2025} introduced a novel neural preconditioner specifically for a class of 3D-1D mixed-dimensional PDEs in which we propose to minimize a residual-based loss functional of the form:
\begin{equation*}
\mathcal{L}_{\text{static}}(\theta) = \frac{1}{|\mathcal{B}|} \mathsf{ \sum_{b^\mu \in \mathcal{B}} \left\| b^\mu - A^\mu \mathcal{N}_\theta( b^\mu, \mu) \right\|}^2,
\end{equation*}
where $\mathcal{B}$ is a suitable collection of vectors, and $\mathcal{N}_\theta$ is a nonlinear operator parameterized by the neural network weights $\theta$. 
This method showcased the potential of unsupervised operator learning via convolutional neural networks (U-Nets~\cite{williams2024unifiedframeworkunetdesign}) to create preconditioners that can adapt to different shapes of a 1D manifold without the need for retraining and can effectively scale to various mesh resolutions.
However, despite its simplicity, this formulation presents limitations, such as limited interpretability in terms of solver dynamics and poor alignment with performance metrics relevant to iterative methods.
To overcome these limitations, we propose an enhanced learning strategy—herein referred to as the \emph{dynamic strategy}—that explicitly incorporates solver convergence dynamics into the training process. The dynamic residual-based loss directly reflects the geometric and convergence properties intrinsic to Krylov subspace methods by employing a differentiable implementation (in the sense of automatic differentiation \cite{gunes2015automatic}) of the Flexible GMRES (FGMRES) algorithm~\cite{saad1993flexible}. Within this framework, the neural preconditioner is influenced by the residual at each iteration. Specifically, we define a novel dynamic loss functional as
\begin{equation*}
\mathcal{L}_{\text{dynamic}}^{(M)}(\theta) = \frac{1}{|\mathcal{B}|} \sum_{\mathsf{b^\mu} \in \mathcal{B}} \sum_{i=1}^{M} \left| s_i(\theta, \mathsf{b^\mu}) \right|,
\end{equation*}
where quantities $s_i(\theta, \mathsf{b^\mu})$, which are closely linked to the algorithm convergence behavior,  offer a geometric characterization of the Krylov subspace generated by the problem $\mathsf{(A^\mu,b^\mu)}$ and preconditioned by the operator $\mathcal{N}_\theta$. The integer $M$ defines the optimization window, focusing the training on the initial iterations of the solver.
Given this construction, we have a structured and solver-aligned loss functional.
In this context, our approach bridges numerical linear algebra with deep learning by leveraging modern automatic differentiation tools (e.g., PyTorch). This framework enables the efficient gradient-based optimization of neural preconditioners through backpropagation across the solver's iterative computations. The combined static-dynamic training strategy enables the neural preconditioner to adapt effectively to complex geometries variations characteristic of mixed-dimensional PDEs, resulting in significantly improved solver performance across diverse problem instances.

The remainder of the paper is organized as follows. Section~2 introduces the mixed-dimensional PDE model and its finite element discretization. Section~3 describes the static pretraining strategy. Section~4 presents foundational concepts and geometric results from Krylov subspace theory that support our dynamic optimization. Sections~5 and 6 develop the dynamic loss formulation and detail the training algorithms using differentiable linear algebra. Sections 7 and 8 report numerical results and conclude with final remarks and future directions.

\section{Problem Setting}

Throughout the manuscript, we adopt a consistent notation scheme to distinguish between different mathematical objects. Calligraphic letters such as $\mathcal{P}$ and $\mathcal{H}$ are used to denote spaces or sets. The notation $|\mathcal{P}|$ refers to the cardinality of the set $\mathcal{P}$. In this context, an exception to this rule is the neural preconditioner $\mathcal{N}_{\theta}$. Finite-dimensional vectors are represented as $\mathsf{v}$ or $\mathsf{w}$, and the same convention applies to matrices, such as $\mathsf{A}$ and $\mathsf{B}$. Given a matrix $\mathsf{A}$, its entry located at row~$i$ and column~$j$ is written as $[\mathsf{A}]_{ij}$. The inner product between two elements is denoted by $\langle \cdot, \cdot \rangle$. We denote vector norms using the Euclidean norm $\lVert \cdot \rVert$.

We consider the numerical solution of a family of parametrized mixed-dimensional elliptic problems representative of a broader class of PDEs defined over coupled domains of different dimensions. Specifically, we study a 3D-1D coupled problem set in a three-dimensional domain $\Omega \subset \mathbb{R}^3$ and a one-dimensional manifold $\Lambda \subset \Omega$, thoroughly analyzed in~\cite{Laurino20192047} and typically modeling slender structures.
The continuous mixed-dimensional formulation reads
\begin{equation}
\label{eq:cont_problem}
\begin{aligned}
- \nabla \cdot (k_\Omega \nabla u_\Omega) + \sigma_\Omega u_\Omega + 2\pi \epsilon \, (\mathcal{T}_\Lambda u_\Omega - u_\Lambda) \, \delta_\Lambda &= 0, && \text{in } \Omega, \\
- \partial_s (k_\Lambda \, \partial_s u_\Lambda) + 2\pi \epsilon \, (u_\Lambda -\mathcal{ T}_\Lambda u_\Omega) &= 0, && \text{on } \Lambda,
\end{aligned}
\end{equation}
subject to appropriate boundary conditions. Here, $u_\Omega$ and $u_\Lambda$ are the unknown fields defined in $\Omega$ and $\Lambda$, respectively; $k_\Omega$, $\sigma_\Omega$, and $k_\Lambda$ are physical coefficients; and $\epsilon > 0$ controls the strength of the coupling between the domains. The operator $\mathcal{T}_\Lambda$ denotes the coupling of $\Omega$ to $\Lambda$, typically implemented as a cross-sectional average in tubular neighborhoods of $\Lambda$. The term $\delta_\Lambda$ indicates that the coupling acts on $\Omega$ as a singular source supported on $\Lambda$.

We discretize \eqref{eq:cont_problem} using the Galerkin projection on a (broken) finite element space $V_h = V_h^\Omega \times V_h^\Lambda$ of dimension $N_h=N_h^\Omega + N_h^\Lambda$, being $N_h^\Omega = \textrm{dim}(V_h^\Omega)$ and $N_h^\Lambda = \textrm{dim}(V_h^\Lambda)$, which are chosen depending on the characteristics of the problem at hand. Such a discretization leads to a family of linear systems of the form
\begin{equation*}
\mathsf{A^\mu u^\mu = b^\mu,} \quad \mu \in \mathcal{P},
\end{equation*}
where $\mu$ indexes the parameter space $\mathcal{P}$ encoding geometric variations of the embedded structure $\Lambda$. The discrete solution vector $\mathsf{u^\mu}$ and right-hand side $\mathsf{b^\mu}$ are partitioned into 3D and 1D components
\begin{equation*}
\mathsf{u^\mu} = \begin{pmatrix} \mathsf{ u_{\mu,0}} \\ \mathsf{u_{\mu,1}} \end{pmatrix}, \quad
\mathsf{b^\mu }= \begin{pmatrix} \mathsf{0} \\ \mathsf{f_{\mu,1}} \end{pmatrix}.
\end{equation*}
The system matrix $\mathsf{A^\mu}$ has the following block structure
\begin{equation*}
\mathsf{A^\mu} = 
\begin{pmatrix} \hspace{4pt}
\mathsf{K_{h,00} + 2\pi \epsilon M_{h,00}^\mu} & \mathsf{2\pi \epsilon M_{h,01}^\mu} \\[7pt]
\mathsf{2\pi \epsilon M_{h,10}^\mu} & \mathsf{K_{h,11}^\mu + 2\pi \epsilon M_{h,11}^\mu }
\hspace{4pt}\end{pmatrix},
\end{equation*}
where the blocks $\mathsf{K_{h,ij}}$ correspond to discretizations of the elliptic operators, and $\mathsf{M_{h,ij}^\mu}$ represent parameter-dependent coupling terms.

The principal computational difficulty lies in solving this system efficiently across multiple instances $\mu \in \mathcal{P}$. In particular, the dominant cost is associated with solving the 3D subproblem
\begin{equation*}
\mathsf{A_{h,00}^\mu u_{\mu,0} = b_{h,0}^\mu },
\end{equation*}
where
\begin{equation*}
\mathsf{A_{h, 00}^\mu := K_{h,00}} + 2\pi \epsilon\mathsf{ M_{h,00}^\mu}, \quad
\mathsf{b_{h,0}^\mu} := -2\pi \epsilon \mathsf{M_{h,01}^\mu u_{\mu,1}}.
\end{equation*}
Note that the relevant right-hand sides for this reduced system belong to
\begin{equation*}
\mathcal{B} := \left\{ -2\pi \epsilon \mathsf{M_{h,01}^\mu x} \;\middle|\; \mu \in \mathcal{P},\; \mathsf{x} \in \mathbb{R}^{N_{\Lambda, h}} \text{ such that } \mathsf{(K_{h,11}^\mu + 2\pi \epsilon M_{h,11}^\mu) x = f_{\mu,1}} \right\}.
\end{equation*}
This set $\mathcal{B}$ characterizes the structure of the right-hand sides that arise during the elimination of the 1D variables from the coupled system and plays a key role in training neural preconditioners.

\section{Static Pretraining for Neural Preconditioning}

In this section, we introduce the static pretraining phase for neural preconditioners, which serves as the initial step in a two-stage learning process for solving families of parametrized linear systems arising from mixed-dimensional PDEs. The primary goal of this phase is to obtain a good initialization for the neural network parameters $\theta$ using an unsupervised learning approach, without requiring knowledge of the true solution vectors. The methodology described here has already been proposed in~\cite{dimola2025}.

From now on, we consider a family of parameter-dependent linear systems of the form
\begin{equation*}
\mathsf{A_{h,00}^\mu u_{\mu,0} = b_{h,0}^\mu },  \quad \mu \in \mathcal{P}\, ;
\end{equation*}
Given the absence of  ambiguity, we will drop all subscripts to ease the notation, leading to,
\begin{equation}
\label{eq:3D_problem}
\mathsf{A^\mu u^{\mu} = b^\mu },  \quad \mu \in \mathcal{P}\, ;
\end{equation}
where $\mathcal{P}$ denotes a compact parameter space encoding geometrical or physical variations of the problem. For each $\mu \in \mathcal{P}$, $\mathsf{A}^\mu \in \mathbb{R}^{N\times N}$ is a symmetric positive-definite (SPD) matrix resulting from the discretization of the mixed-dimensional PDEs, and $\mathsf{b^\mu} \in \mathbb{R}^N$ is the corresponding right-hand side.
Notably, although the problem at hand involves an SPD matrix, we do not exploit this property in the solver for two main reasons: $i)$ Eq.~\eqref{eq:3D_problem} represents a specific instance of a broader class of operators that may lack symmetry, such as those arising in transport-dominated regimes; $ii)$ learned methods like our neural preconditioner generally break symmetry unless such structure is explicitly enforced in the architecture (however, this is a capability for which no efficient implementation is currently available).

The neural preconditioner is defined as a nonlinear operator $\mathcal{N}_{\theta}: \mathbb{R}^N \times \mathcal{P} \rightarrow \mathbb{R}^N$, represented by a neural network with parameters $\theta$. The aim is to train $\mathcal{N}_{\theta}$ so that it approximately inverts the operator $\mathsf{A^\mu}$ over a subset of right-hand sides relevant to the problem.
To this end, we minimize the following residual-based, unsupervised loss functional
\begin{equation*}
\mathcal{L}_{\text{static}}(\theta) = \frac{1}{|\mathcal{P}|} \sum_{j=1}^{|\mathcal{P}|} \frac{1}{|\mathcal{K}_{\mu_j}|} \sum_{v \in \mathcal{K}_{\mu_j}} \left\| \mathsf{v} - \mathsf{A^{\mu_j}} \mathcal{N}_{\theta}(\mathsf{v}, \mu_j) \right\|^2,
\end{equation*}
where $\{\mu_j\}_{j=1}^{|\mathcal{P}|} \subset \mathcal{P}$ are sampled parameters, and $\mathcal{K}_{\mu_j} \subset \mathbb{R}^N$ is a training set of normalized right-hand sides associated with parameter $\mu_j$.
Each training set $\mathcal{K}_{\mu_j}$ includes both physics-informed and spectrally enriched right-hand sides, and is defined as
\begin{equation*}
\mathcal{K}_{\mu_j} = \left\{\mathsf{ \frac{b^{\mu_j}}{\|b^{\mu_j}\|}} \right\} \cup \mathcal{D}_{\mu_j},
\end{equation*}
where $\mathsf{b_{\mu_j}}$ is a canonical forcing term associated with $\mu_j$, and $\mathcal{D}_{\mu_j}$ is a data augmentation set constructed to enrich the spectral content of the training data.\\
To construct the set $\mathcal{D}_{\mu_j}$, we sample unrelated random unit vectors $\mathsf{r_{k,h}^{\mu_j}} \in \mathbb{S}^{N_h-1} \subset \mathbb{R}^{N_h}$ uniformly from the unit hypersphere. This is achieved in practice by drawing a Gaussian random vector $\mathsf{v} \sim \mathcal{N}(\mathbf{0}, \mathbf{I}_{N_h})$ and projecting it onto the sphere
\begin{equation*}
\mathsf{r^{\mu_j}_{k,h}} = \frac{\mathsf{v}}{\|\mathsf{v}\|}\, .
\end{equation*}
This procedure generates samples that are uniformly distributed on $\mathbb{S}^{N_h-1}$ due to the rotational invariance of the standard multivariate normal distribution; since the distribution of $\mathsf{v}$ is isotropic, the normalization $\mathsf{v}/\|\mathsf{v}\|$ yields a uniform distribution over the unit sphere.
As such, the resulting augmented set is given by
\begin{equation*}
\mathcal{D}_{\mu_j} = \left\{ \mathsf{r_{0,h}^{\mu_j}},\, \mathsf{r_{1,h}^{\mu_j}},\, \ldots \right\}.
\end{equation*}
This spectral enrichment ensures that the neural preconditioner learns to handle both low- and high-frequency components of the operator spectrum, which are essential for effective preconditioning in mixed-dimensional settings.
The minimization of the empirical risk $\mathcal{L}_{\text{static}}$ is performed using standard gradient-based optimization methods such as Adam. Importantly, this static pretraining phase is fully unsupervised—it does not rely on access to exact solutions $\mathsf{u^{\mu}}$, thereby significantly reducing the computational cost associated with the offline training.
The resulting neural preconditioner $\mathcal{N}_{\theta}$ provides a robust initial approximation of the inverse action of $\mathsf{A}^\mu$ on representative right-hand sides. This initialization is then refined during a subsequent solver-integrated fine-tuning stage, where the learned operator is further adapted to enhance the convergence of Krylov subspace methods such as GMRES.

\section{Geometric Aspects of Krylov Subspace Methods}
This section outlines the geometric foundations of Krylov subspace methods, following the rigorous treatment in~\cite{Eiermann_Ernst_2001}, with emphasis on the convergence-related quantities that motivate our dynamic loss formulation for neural preconditioning.

We consider the linear system:
\begin{equation*}
\mathsf{ A u = f, \quad A} \in \mathbb{R}^{n \times n}, \quad \mathsf{f} \in \mathbb{R}^n\, ,
\end{equation*}
with an initial guess $\mathsf{u_0}$ and residual $\mathsf{r_0 = f - A u_0}$. Krylov subspace methods generate iterates $\mathsf{u_j \in u_0 +} \mathcal{C}_j$, where the \textit{correction subspace} is defined by,
\begin{equation*}
\mathcal{C}_j = \mathcal{K}_j(\mathsf{A, r_0}) := \text{span}\{\mathsf{r_0, A r_0, \dots, A^{j-1} r_0}\}.
\end{equation*}
Letting $\mathcal{W}_j:=\{\mathsf{x}\,  |\,\mathsf{x=Ay, \,y}\in \mathcal{C}_j \}$ denote the image subspace of the correction subspace, with respect to matrix $\mathsf{A}$, we get that in minimal residual Krylov subspace methods, the residual at iteration $j$ is enforced to be orthogonal to $\mathcal{W}_j$, i.e.,
\begin{equation*}
\mathsf{r_j = f - A u_j} \perp \mathcal{W}_j\, .
\end{equation*}
The iterate $\mathsf{u_j = C_j \nu^*_j}$ is thus obtained by solving the minimization problem:
\begin{equation*}
\mathsf{\nu^*_j = \arg\min_{\nu \in \mathbb{R}^j} \|r_0 - A C_j} \nu\|\, ,
\end{equation*}
where $\mathsf{C_j} \in \mathbb{R}^{n\times j}$ is a column basis matrix for correction space $\mathcal{C}_j$.

As expected, the convergence of Krylov methods is intimately connected to the geometry of residual vectors relative to Krylov subspaces. Indeed, given a nonzero vector $\mathsf{x}$ and a subspace $\mathcal{Y}$, define the the principal angle between them by
\begin{equation*}
\angle(\mathsf{x}, \mathcal{Y}) := \arcsin \left( \frac{\|(I - P_{\mathcal{Y}})\mathsf{x}\|}{\|\mathsf{x}\|} \right)\, ,
\end{equation*}
where $P_{\mathcal{Y}}$ denotes the orthogonal projector onto $\mathcal{Y}$;  then $s_j :=\sin (\angle(\mathsf{r_{j-1}}, \mathcal{W}_j))$ denote the sine of the angle between $\mathsf{r_{j-1}}$ and $\mathcal{W}_j$.
Given that, by construction, the subspace sequence $\{\mathcal{W}_j\}$ is nested:
\begin{equation*}
\{0\}\equiv\mathcal{W}_0 \subset \mathcal{W}_1 \subset \ldots \subset \mathcal{W}_{j-1} \subset \mathcal{W}_j\, .
\end{equation*}
It can be proved that the residual norm satisfies the following recurrence,
\begin{equation}
\|\mathsf{r_j}\| = s_j \|\mathsf{r_{j-1}}\| = s_j s_{j-1} \cdots s_1 \|\mathsf{r_0}\|\, ,
\label{eq:sinus_recurrence}
\end{equation}
so that, fast convergence \textit{corresponds} to small angle sines $|s_j| \ll 1$ at early iterations.

The key computational tool in Krylov methods is the Arnoldi process, which constructs an orthonormal basis $\mathsf{V_j}$ of $\mathcal{K}_j(\mathsf{A, r_0})$ that satisfies the matrix relation
\begin{equation*}
\mathsf{A V_j = V_{j+1} \widetilde{H}_j}\, ,
\end{equation*}
with $\mathsf{\widetilde{H}_j }\in \mathbb{R}^{(j+1) \times j}$ an upper Hessenberg matrix. Geometric information is hidden in the Hessenberg matrix, so that, given an orthogonal matrix $\mathsf{Q_j} \in \mathbb{R}^{j+1, j+1}$  such that
\begin{equation*}
\mathsf{Q_j \widetilde{H}_j} = 
\begin{bmatrix} 
\mathsf{R_j} \\ 
\mathsf{0} 
\end{bmatrix}.
\end{equation*}

We have that 
\begin{equation*}
s_m=\left|\frac{\mathsf{[Q_j]}_{ j+1,1}}{[\mathsf{ Q_{j-1}}]_{j,1}}\right|.
\end{equation*}
If one constructs the orthogonal transformation as a sequence of Givens rotations $\mathsf{Q_j = G_j G_{j-1} \cdots G_1}$, where each $\mathsf{G_k}$ annihilates subdiagonal elements of $\mathsf{\widetilde{H}_j}$, we can recover $s_j$ in $\mathsf{G_j}$ skew-symmetric part.

Given all stated above, consider a $\theta$-parametrized, possibly non-linear, preconditioner $\mathcal{N_\theta}: \mathbb{R}^n \rightarrow \mathbb{R}^n$ and the associated flexible variant GMRES, so that the associated Krylov subspace generated in the iterative process results 
\begin{equation*}
\mathcal{K}_j^{\text{FGMRES}} = \text{span}\{\mathsf{r}_0, \mathsf{A} \circ  \mathcal{N}_{\theta}(\mathsf{r}_0), \dots, (\mathsf{A} \circ \mathcal{N}_{\theta})^{j-1}(\mathsf{r}_0)\}\, ,
\end{equation*}
and the Arnoldi relation generalizes to
\begin{equation*}
\mathsf{A Z_j = V_{j+1} \widetilde{H}_j}\, ,
\end{equation*}
where $\mathsf{Z_j} = [\mathcal{N}_{\theta}(\mathsf{v}_0), \dots, \mathcal{N}_{\theta}(\mathsf{ v_{j-1}})]$, with $\{\mathsf{v_l\}_{l=1}^j}$ the orthonormal basis vectors generated during the iterations. 
 As in the classical case, the sine angles $s_j$ can still be computed from $\mathsf{\widetilde{H}_j}$ via Givens rotations, but, given the introduction of $\mathcal{N_\theta}$ in the nested subspace construction process, they become differentiable functions of the operator parameters $\theta$, thus we write $s_{j,\theta}$.

These geometric insights form the basis of the dynamic loss functional employed in our neural preconditioning framework. By searching for optimal parameters $\theta^\star$ penalizing large values of $\{|s_{j,\theta}|\}$, the dynamic loss explicitly minimizes the misalignment between residual vectors and Krylov search directions, enabling solver-aware optimization aligned with the convergence behavior of FGMRES.

\section{Dynamic Loss Functional Formulation}

We now formalize the geometric intuition developed in the previous section into a differentiable optimization framework suitable for training $\mu$-parametrized nonlinear preconditioners via gradient-based methods. Our approach is centered on the construction of a loss functional that quantifies and penalizes the misalignment between the residual vector and the Krylov subspace generated during the solution process, to accelerate convergence.
Such a misalignment is measured using the sine of the angle between the residual and the Krylov subspace at each iteration. This metric is embedded directly into the training procedure and is evaluated dynamically during the execution of the Flexible GMRES (FGMRES) algorithm, wherein each Krylov basis vector is modified by a nonlinear, parameterized preconditioner $\mathcal{N}_{\theta}: \mathbb{R}^N \times \mathcal{P} \rightarrow \mathbb{R}^N$.
 For the parameter instance, $\mu$, consider $M$ step of the generalised Arnoldi procedure,
\begin{equation*}
    \mathsf{A^\mu Z_M^\mu = V_{M+1}^\mu \widetilde{H}_M^\mu}\, ,
\end{equation*}
with
\begin{equation*}
   \mathsf{Z_M^\mu} = [\mathcal{N}_{\theta}(\mathsf{v}_0, \mu),\, \dots,\, \mathcal{N}_{\theta}(\mathsf{v_{j}}, \mu),\, \ldots,\, \mathcal{N}_{\theta}(\mathsf{v_{M-1}}, \mu)].
\end{equation*}
Each preconditioned vector $\mathsf {A^\mu}\mathcal{N}_\theta({\mathsf{v}_j}, \mu)$ undergoes a Gram–Schmidt orthogonalization process, preserving the nested structure of the Krylov subspaces.

At iteration $M$, we define $\{s^\mu_1, \ldots, s^\mu_M\}$ the sequence of sines of the angles between residuals and Krylov subspaces generated up to that point.
This sequence is a by-product of the sequence of $M$ Givens rotations, $\mathsf{G_M G_{M-1}\ldots G_2 G_1}$ that orthogonalize $\mathsf{\widetilde{H}_M^\mu}$. We envelop this process of sine extraction from the Givens rotations orthogonalizing sequence with the shorthand notation $\textbf{AG}_M$. 

Since the basis vectors depend on $\theta$ via the action of $\mathcal{N}_\theta$, the entire sequence $\{s_j^\mu\}_{j=1}^M \equiv \{s_{j, \theta}^\mu\}_{j=1}^M$ becomes differentiable with respect to the model parameters.
The goal is to find a parameter vector $\theta^\star$ such that
\begin{equation*}
    \prod_{j=1}^{M} |s^\mu_{j,\theta^\star}|\ll 1\, ,
\end{equation*}
for all $\mu \in \mathcal{P}$, indicating rapid convergence of the residual norm over the first $M$ iterations for the whole  problem family.

\begin{figure}
    \centering
    \includegraphics[width=0.9\linewidth]{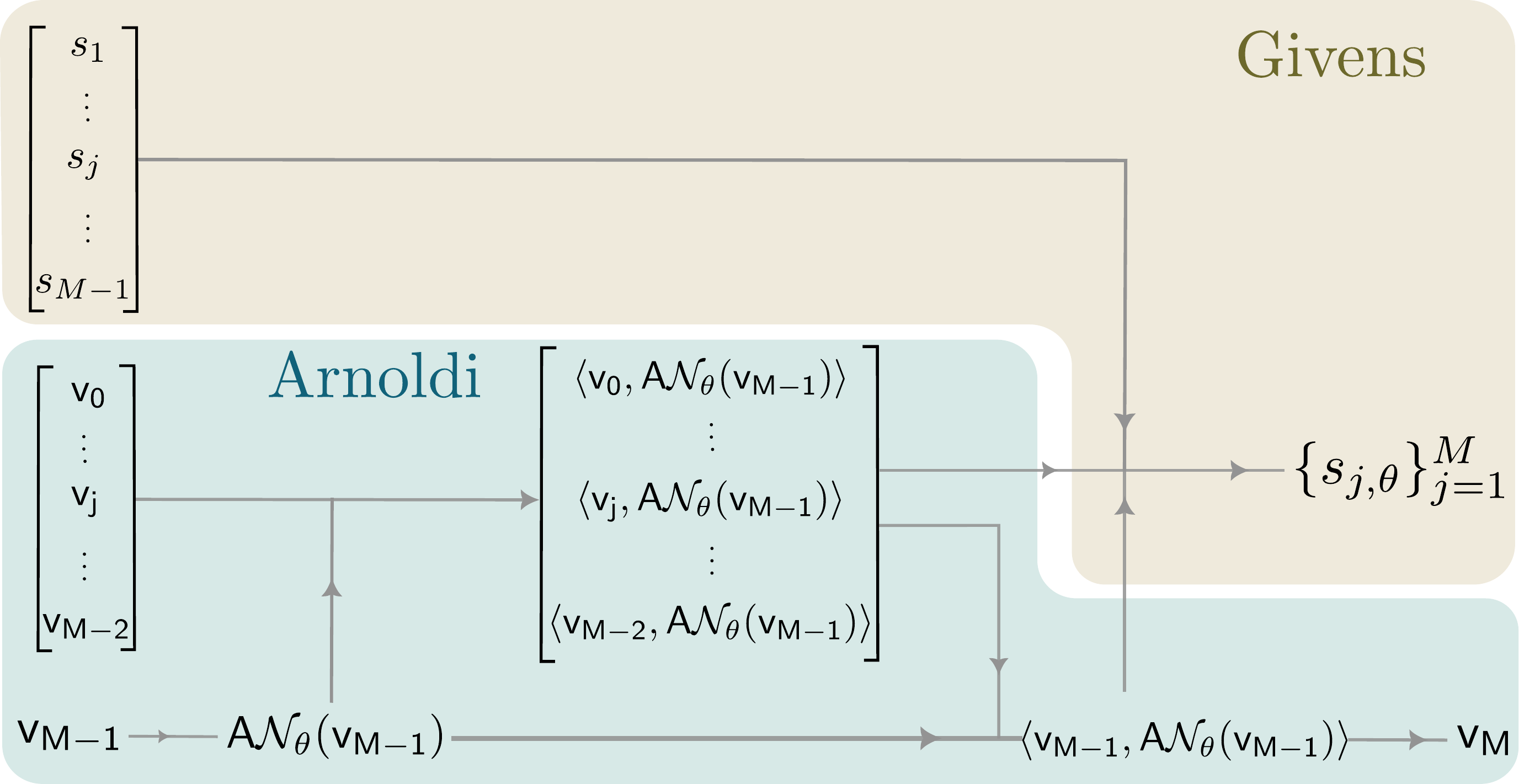} \vspace{10pt}
    \caption{Computational graph corresponding to the $M$-th iteration of the Arnoldi-Givens process ($\textbf{AG}_M$). Dependencies on the learnable parameters $\theta$ are introduced through the nonlinear preconditioner $\mathcal{N}_\theta$. For brevity, we suppress the dependence on parameters $\mu$ and $\theta$ when superfluous. }
    \label{fig:DAG}
\end{figure}

The recursive nature of the Arnoldi-Givens process induces a computational graph where each variable depends on the previous ones. This structure forms a Directed Acyclic Graph (DAG), in which nodes represent intermediate quantities (e.g., inner products, basis vectors, rotation angles), and edges encode standard algebraic operations, such as those involved in performing orthogonalization or applying Givens rotations. The parameter dependence enters the graph via the preconditioner $\mathcal{N}_\theta$, which affects every step of the $\textbf{AG}_M$ iterative procedure. In Fig.~\ref{fig:DAG}, we visualize this computational structure, emphasizing the role of $\theta$ in shaping the algorithm’s trajectory.

This computational structure enables seamless application of reverse-mode automatic differentiation. In particular, frameworks such as PyTorch automatically construct and traverse the computational graph, allowing efficient evaluation of the gradient \cite{paszke2019pytorch}
\begin{equation*}
    \delta \theta \mapsto \left\{ \partial s^\mu_{j,\theta}/\partial \theta \right\}_{j=1}^M\, ,
\end{equation*}
which can be used to train the neural preconditioner by minimizing a suitably defined loss functional.
The complete learning procedure is outlined in Algorithm 1.
\begin{algorithm}
\label{alg:dynamic}
\small
\caption{Training a Neural Preconditioner via Krylov Subspace Alignment}
\begin{algorithmic}[1]
\State Select linear systems $(\mathsf{ A^\mu,b^\mu)}$, for $\mu \in \mathcal{P_{\text{train}}}$
\State Initialize the neural operator $\mathcal{N}_\theta$ and the optimizer
\For{each training epoch}
    \State Compute the initial residual $\mathsf{r_0^\mu}$
    \For{$j = 1$ to $M$}
        \State Construct $ \mathsf{A^\mu Z_j^\mu = V_{j+1}^\mu \widetilde{H}_j^\mu}\, $ \hspace{105pt}\rlap{\smash{$\left.\begin{array}{@{}c@{}}\\{}\\{}\\{}\end{array}\color{black}\right\}
          \color{black}\begin{tabular}{l}\small $\textbf{AG}_M(\mathsf{A^\mu, r_0^\mu}, \mathcal{N}_\theta)$.\end{tabular}$}}
        \State Orthogonalize $\mathsf{\widetilde{H}_j^\mu}$ via Givens rotations; compute $s^\mu_{j,\theta}$  
    \EndFor
    \State Evaluate the loss $\mathcal{L}(\theta) = \sum_{j=1}^M \left\langle |s^\mu_{j,\theta}| \right\rangle_\mu$
    \State Compute $\nabla_\theta \mathcal{L}$ using backpropagation
    \State Update $\theta \leftarrow \theta - \eta \nabla_\theta \mathcal{L}$
\EndFor

\end{algorithmic}
\end{algorithm} 
The average $\left\langle \cdot \right\rangle_\mu$ denotes the empirical mean over a batch of parameterized linear systems $(\mathsf{A^\mu}, \mathsf{b^\mu})\, , \mu \in \mathcal{P}_{\text{train}} \subset \mathcal{P}$ or $\mu \in \mathcal{P}_\text{test}\subset \mathcal{P}$, depending on the context,
\begin{equation*}
    \left\langle \mathsf{x^\mu} \right\rangle_\mu := \frac{1}{ |\mathcal{P}|} \sum_{\mu \in \mathcal{P}} \mathsf{x^\mu}\, .
\end{equation*}
This formulation allows the learned preconditioner to generalize across the entire problem class.
The dynamic loss functional is thus defined as \vspace{-5pt}
\begin{equation*}
    \mathcal{L}_{\text{dynamic}(M)}(\theta) := \frac{1}{M}\sum_{j=1}^M \left\langle |s^\mu_{j,\theta}| \right\rangle_{\mu \in \mathcal{P}_\text{train}},\vspace{-7pt}
\end{equation*}
where the sine series $\{s^\mu_{j,\theta}\}$ is evaluated \textit{dynamically} in the training process, i.e. whenever parameter is updated, $\theta \leftarrow (\theta-\eta \nabla_\theta\mathcal{L})$ (line 11 in Algorithm \ref{alg:dynamic}), the Arnoldi-Givens algorithm is called to obtain the new series 
\begin{equation*}
    \textbf{AG}_M: (\mathsf{A^\mu, r_0^\mu},\mathcal{N}_{(\theta- \eta \nabla_\theta\mathcal{L})} ) \longmapsto \{s_{j,(\theta- \eta \nabla_\theta\mathcal{L})}^\mu\}_{j=1}^M.    
\end{equation*}
As $\mathcal{L}_{\text{dynamic}(M)} \approx \mathcal{O}(M)$, it provides a robust surrogate for minimizing residual misalignment over the first $M$ iterations.
Moreover, the inequality,
\begin{equation*}
    \prod_{j=1}^M \langle |s^\mu_j| \rangle_\mu \leq \left(  \mathcal{L}_{\text{dynamic}(M)}(\theta) \right)^M\, ,
\end{equation*}
relates the geometric mean of the sine values to the average loss, offering a theoretical guarantee that minimizing $\mathcal{L}_{\text{dynamic}(M)}$ leads to improved convergence behavior.

The parameter $M$ plays a central role in balancing optimization depth and computational efficiency. While larger values of $M$ provide a more accurate reflection of long-term solver dynamics, they also result in deeper computational graphs and increased memory usage. In practice, moderate values of $M$ (e.g., $5 \leq M \leq 10$) are sufficient to achieve meaningful learning without incurring prohibitive computational costs.
It is important to note that, unlike matrix inversion, preconditioning does not require a unique operator. The space of valid preconditioners is large and diverse, allowing the learning process to exploit non-uniqueness to improve convergence.
In the limiting case $M = 1$, the loss effectively trains $\mathcal{N}_\theta$ to approximate the inverse of $\mathsf{A}^\mu$, which is typically insufficient for practical tolerances. Larger values of $M$ enable iterative refinement and allow the network to learn effective strategies for accelerating convergence across a wide range of systems.

\section{Training Strategy for the Neural Preconditioner}
\label{sec:training}

The training of the nonlinear preconditioner $\mathcal{N}_\theta$ follows a two-phase procedure designed to balance computational efficiency with solver-aware optimization. First, we perform an unsupervised pretraining stage based on a static loss functional that avoids unrolling the Krylov solver. Then, the network undergoes fine-tuning using a dynamic loss that incorporates geometric information derived from the iterative process.
This two-phase strategy reflects a trade-off between generalization and specialization. The static phase enables the network to learn a coarse but broadly applicable approximation of the inverse operator, leveraging only residual information without requiring ground-truth solutions. Data augmentation-based randomized perturbations enrich the spectral content of the training set, promoting robustness across parameter variations.
On the other hand, the subsequent dynamic phase injects solver-specific information into the training loop by minimizing the absolute value of sine angles $s_{j, \theta}$ between residuals and Krylov subspaces. These angles govern convergence behavior in GMRES and its flexible variant, and their minimization aligns the learned preconditioner with the subspace geometry of the iterative process. This fine-tuning step optimizes the early-stage convergence rate and reduces iteration counts in practice. By decoupling representation learning from solver dynamics, this hybrid approach ensures both structural adaptability and performance-driven refinement. The static loss captures coarse inverse structure, while the dynamic loss promotes alignment with Krylov geometry, resulting in a preconditioner that is both generalizable and solver-aware.

\subsection{Neural Architecture: A Multi-Level U-Net Design}

We briefly recall here the definition of the architecture chosen for $\mathcal{N}_\theta$, previously discussed in~\cite{dimola2025}. To define the hypothesis space $\mathcal{H}$, we adopt a neural operator $\mathcal{N}_\theta$ instantiated via a U-Net architecture of depth $L$, denoted by $U_L$. This architecture is well-suited for mixed-dimensional PDEs due to its ability to capture multiscale spatial features.
Then, we set $\mathcal{H} = \{ U_L(\cdot\,;\theta) \}_{\theta \in \mathbb{R}^t; \, L \in \mathbb{N}^+}$, so that the neural preconditioner is defined as $\mathcal{N}_\theta \equiv U_L(\cdot\,;\theta)$.

The U-Net architecture comprises an encoder,
\begin{equation*}
\Phi = \Phi_0 \circ \Phi_1 \circ \dots \circ \Phi_{L-1}: \mathbb{R}^{c_{\text{in}} \times n} \to \mathbb{R}^{c_b \times m}\, ,
\end{equation*}
and a decoder,
\begin{equation*}
\Psi = \Psi_{L-1} \circ \dots \circ \Psi_0: \mathbb{R}^{c_b \times m} \to \mathbb{R}^{c_{\text{out}} \times n}.
\end{equation*}
Given an input tensor $\mathbf{X} = [\mathbf{X}_1 \mid \ldots \mid \mathbf{X}_{c_{\text{in}}}]$ with $\mathbf{X}_i \in \mathbb{R}^{n_1 \times \dots \times n_d}$, the recursive structure of the U-Net is defined as
\begin{align*}
U_0(\mathbf{X}) &:= \Psi_0 \circ \Phi_0(\mathbf{X}), \\
U_j(\mathbf{X}) &:= \Psi_j\left( \left[ U_{j-1} \circ \Phi_j(\mathbf{X}) \;\middle|\; \Phi_j(\mathbf{X}) \right] \right), \quad j > 0\, ,
\end{align*}
where $[\;\cdot \mid \cdot\;]$ denotes channel-wise concatenation (skip connection), enhancing both gradient propagation and spatial information retention.

Each convolutional block $\Phi_j$ and $\Psi_j$ applies $c_{\text{out}}$ filters ${\bf k}_{i,j}$ to the $c_{\text{in}}$ input channels,
\begin{equation*}
\mathbf{Y}_i = \sum_{j=1}^{c_{\text{in}}} {\bf k}_{i,j} \ast \mathbf{X}_j, \qquad i = 1, \dots, c_{\text{out}}\, ,
\end{equation*}
with a discrete convolution defined by
\begin{equation*}
({\bf k}_{i,j} \ast \mathbf{X}_j)(\mathbf{p}) = \sum_{\mathbf{q} \in \mathbb{Z}^d} {\bf k}_{i,j}(\mathbf{q})\, \mathbf{X}_j(\mathbf{p} - \mathbf{q}).\vspace{-7pt}
\end{equation*}

To comply with convolutional formats, the raw inputs $\{\mathbf{x}_i\}_{i=1}^{c_{\text{in}}} \subset \mathbb{R}^{N_h}$ are reshaped into a tensor
\begin{equation*}
\mathbf{x}_1, \ldots, \mathbf{x}_{c_{\text{in}}} \mapsto \mathbf{X} \in \mathbb{R}^{c_{\text{in}} \times n_1 \times \dots \times n_d}, \qquad N_h = n_1 \cdot \dots \cdot n_d.
\end{equation*}
This choice limits the architecture to structured tensor-product meshes; generalization to unstructured domains may require mesh-aware models~\cite{FrancoMINN}.

The specific architecture used in our experiments is detailed in Table~\ref{tab:unet}. It includes downsampling (via max pooling) and upsampling (via transposed convolution) layers.
\begin{table}[h!]
    \centering
    \begin{adjustbox}{width=\textwidth}
    \footnotesize
    \begin{tabular}{lccccc}
    \hline
        \textbf{Layer (type)} & \textbf{Input size} & \textbf{Output size} & \textbf{Input channels} & \textbf{Output channels} & \textbf{\# Parameters} \\ \hline \\
        Conv. & $21^3$ & $21^3$ & 2 & 32 & 13,840 \\ 
        Conv. + Max Pooling & $21^3$ & $10^3$ & 32 & 64 & 55,296 \\ 
        Conv. + Max Pooling & $10^3$ & $5^3$ & 64 & 128 & 221,184 \\ 
        Conv. & $5^3$ & $5^3$ & 128 & 128 & 442,368 \\ 
        Transposed Conv. & $5^3$ & $10^3$ & 128 & 64 & 286,784 \\ 
        Transposed Conv. & $10^3$ & $21^3$ & 64 & 32 & 71,912 \\ \\ \hline
    \end{tabular}
    \end{adjustbox}
    \caption{ Architecture of the three-level U-Net $U_3$. The total number of trainable parameters is approximately $10^6$.}
    \label{tab:unet}
\end{table}

\subsection{Phase 1: Static Pretraining via Residual Loss}
\label{sec:static_training}
In the static pretraining phase, we minimize a residual-based loss functional that does not require unrolling the Krylov solver
\begin{equation*}
\mathcal{L}_{\text{static}}(\theta) = \frac{1}{|\mathcal{P}_\text{train}|} \sum_{j=1}^{|\mathcal{P}_\text{train}|} \frac{1}{|\mathcal{K}_{\mu_j}|} \sum_{\mathsf{v} \in \mathcal{K}_{\mu_j}} \left\| \mathsf{v} - \alpha^{-1} \mathsf{A}^{\mu_j} \mathcal{N}_\theta(\mathsf{v}, \mu_j) \right\|^2,
\end{equation*}
Here, $\alpha$ is a normalization factor based on the average RMS of the matrix diagonals
\begin{equation*}
\alpha := \left(\frac{1}{|\mathcal{P}_\text{train}|} \sum_{j=1}^{|\mathcal{P}_\text{train}|} \frac{1}{\sqrt{N_h}} \left\| \text{diag}(\mathsf{A}^{\mu_j}) \right\|_2 \right)^{-1}.
\end{equation*}

Inputs consist of the right-hand side vector and the discrete parameter field $\mathsf{d}^{\mu_j}$, reshaped to match the U-Net structure.  The training configuration is summarized in Table~\ref{table:training_data}. We will refer to the neural preconditioner resulting from minimization of $\mathcal{L}_\text{static}$ as $\mathcal{N}_\theta^\text{st}$.  

\begin{table}[h!]
    \centering
    \footnotesize
    \begin{tabular}{l c c | l c c }
    \hline
        \textbf{Training} & \textbf{1D Graphs} & \textbf{Samples} & \textbf{Validation} & \textbf{1D Graphs} & \textbf{Samples} \\ \hline
        Size & 480 & 2400 & Size & 120 & 600 \\
        Batch Size &  & 5 & Epochs &  & 250 \\
        Learning Rate &  & $10^{-3}$ & Scheduler &  & Decaying \\ \hline
    \end{tabular}
    \caption{ Pretraining dataset and parameters. Data augmentation increases the number of samples by a factor of five.}
    \label{table:training_data}
\end{table} \vspace{-30pt}

\subsection{Phase 2: Dynamic Fine-Tuning via Krylov Geometry}

To integrate solver-specific information, we fine-tune the pre-trained model using the dynamic loss functional aligned with the geometry of Krylov  subspaces \vspace{-7pt}
\begin{equation}
\mathcal{L}_{\text{dynamic}(M)}(\theta) = \frac{1}{M}\sum_{j=1}^M \left\langle |s^\mu_{j,\theta}| \right\rangle_{\mu\in \mathcal{P}_\text{train}},  \vspace{-7pt}
\label{eq:dynamic_loss}
\end{equation} 
where the sequence $\{s^\mu_{j,\theta}\}$ of the angles sine between residuals and Krylov subspaces comes from the on-the-fly application of $\textbf{AG}_M$ algorithm and the average $\langle \cdot \rangle_\mu$ is computed over mini-batches of systems $(\mathsf{A^\mu, b^\mu)}$. We use $M = 10$ to prioritize early-stage convergence. The fine-tuning configuration is reported in Table~\ref{table:finetuning_data}. We will refer to the neural preconditioner resulting from minimization of $\mathcal{L}_\text{dynamic}$ as $\mathcal{N}_\theta^\text{dy}$.  
\begin{table}[h!]
    \centering
    \footnotesize
    \begin{tabular}{l@{\hskip 40pt}  c@{\hskip 40pt} | l@{\hskip 40pt}  c@{\hskip 40pt} }
    \hline
        \multicolumn{2}{c}{\textbf{Training}}  & \multicolumn{2}{c}{\textbf{Validation}} \\ \hline
        1D Graphs & 100 & 1D Graphs & 20 \\
        Mini-Batch Size & 20 & Epochs & 150 \\
        Learning Rate & $10^{-3}$ & $M$ (iterations) & 10 \\ \hline
    \end{tabular}
    \caption{\small Fine-tuning dataset and parameters. No data augmentation is used in this phase.}
    \label{table:finetuning_data}
\end{table}

\section{Results and Discussion}

We present a comprehensive analysis of the proposed neural preconditioning strategy. We begin by evaluating the impact of the static pretraining phase on residual reduction. Then, we assess the fine-tuning stage and show how the Krylov geometry-aware loss improves early-stage convergence. Finally, we compare the FGMRES iteration counts with and without preconditioning across a variety of test parameter instances.

\subsection{Pre-Training Step: Static Loss}

We first examine the performance of the neural preconditioner trained with the static loss functional $\mathcal{L}_{\text{static}}$, as introduced in~\cite{dimola2025}. Training was conducted over 100 epochs on a dataset enriched with high-frequency perturbations to promote spectral diversity.

\begin{figure}[h!]
    \centering
    \includegraphics[width=0.99\linewidth]{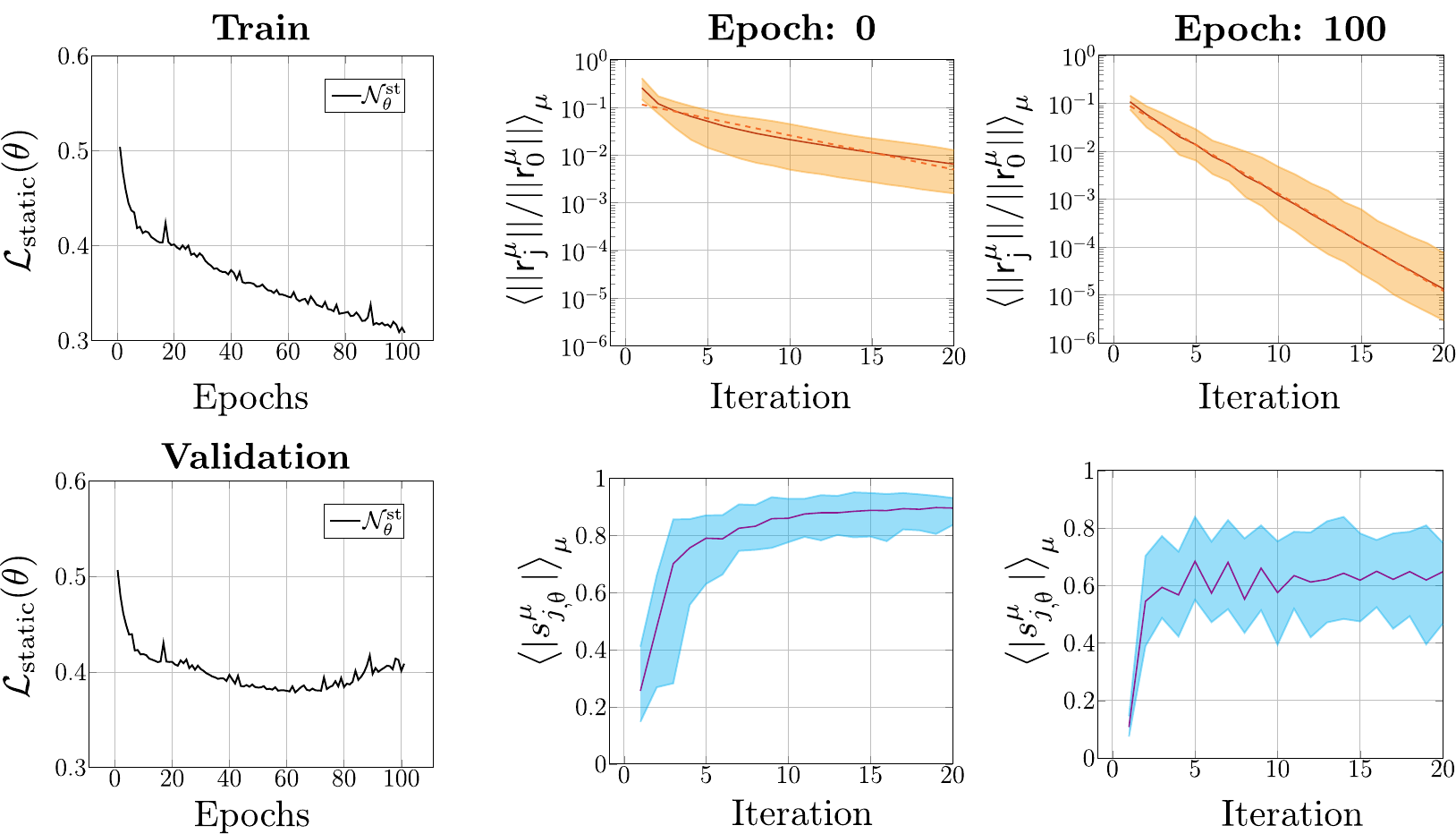}
    \caption{ \textbf{Left column:} Training (top) and validation (bottom) errors obtained during training of $\mathcal{N}_\theta^\text{st}$ preconditioner with residual-based static loss. Training is stopped at the onset of data overfitting. \\\textbf{Middle columns:} Relative residual $\mathsf{||r_j||/||r_0||}$  and averaged absolute subspaces angles sine $\langle|s^\mu_{j,\theta}|\rangle_\mu$ for iteration $j \in  [1, 20]$ at the onset of training (epoch 0). \textbf{Right column:} Relative residual $\mathsf{||r_j||/||r_0||}$  and averaged absolute subspace angles sine $\langle|s^\mu_{j,\theta}|\rangle_\mu$ for iteration $j \in  [1, 20]$ at the end of training step 1 (epoch 100).}
    \label{fig:step1_test}
\end{figure}

As illustrated in Figure~\ref{fig:step1_test} (Left column), both training and validation losses consistently decrease up to approximately 60 epochs, after which a mild overfitting trend emerges. The validation loss stabilizes around a value of 0.4, which is sufficient to produce an effective preconditioner for the considered 3D–1D problem. 

Further analysis highlights the improvement in problem conditioning induced by the learned preconditioner. Before training, the relative residual exhibits slow decay, with sine values close to $0.9$, indicating near-orthogonality of the residual and Krylov image subspaces—an unfavorable condition for convergence (Figure~\ref{fig:step1_test}, Middle column). 

After 100 optimization steps, the average sine value is reduced to approximately $0.6$, resulting in a significantly enhanced convergence profile (Figure~\ref{fig:step1_test}, Right column). As a result, the solver reaches the target tolerance of $10^{-6}$ in roughly 26 iterations, compared to 147 iterations in the unpreconditioned case.

\subsection{Fine-tuning via Krylov Geometry}

\begin{figure}[h]
    \centering
    \includegraphics[width=0.8\linewidth]{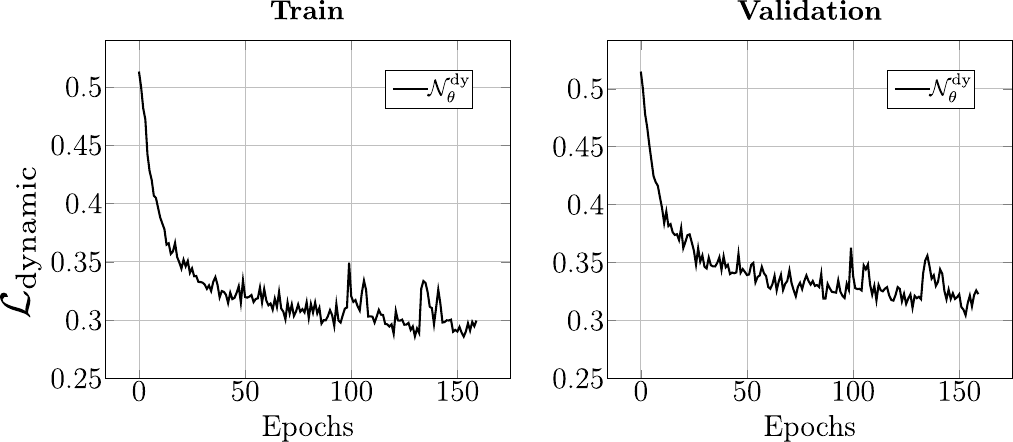}
   
    \caption{Training  and validation loss $\mathcal{L}_\text{dynamic}= \frac{1}{M}\sum_{j=1}^M \langle |s^\mu_{j, \theta}| \rangle_\mu $ during fine-tuning stage for, respectively, $\mu \in \mathcal{P_\text{train}}$ and $\mu \in \mathcal{P_\text{validation}}$.  }
    \label{fig:step2_loss}
\end{figure}

The second phase involves fine-tuning the network using the dynamic loss $\mathcal{L}_{\text{dynamic}}$, which incorporates geometric information extracted from the Krylov solver. As shown in Figure~\ref{fig:step2_loss}, both training and validation losses decrease and stabilize around a value of approximately 0.3, indicating consistent learning behavior across the parameter space and successful integration of subspace alignment objectives.

\subsubsection{Evolution of Principal Angles}
\label{sec:angles}

To gain further insight, we monitor the evolution of the principal angles between the residual vector and the generated Krylov subspace basis vectors, as encoded by the sine quantities $\{|s_{j,\theta}^\mu|\}$. In the ideal case of rapid convergence, these angles should decrease rapidly toward zero, indicating that new Krylov directions contribute substantial progress toward the solution subspace.
Under the first-stage training, based on the static residual loss $\mathcal{L}_{\text{static}}$, the learned preconditioner produces a global improvement of the system conditioning, reflected in a near-uniform reduction of the angles $\{|s_{j,\theta}^\mu|\}$. This behavior suggests that the static training implicitly regularizes the spectrum, albeit without explicit control over the subspace geometry.
In contrast, the second-stage fine-tuning, based on the dynamic loss $\mathcal{L}_{\text{dynamic}}$, directly penalizes large angles within the first $M=10$ Arnoldi iterations.  

\begin{table}[t]
\label{tab:sinus}
\caption{Evolution of principal angles over training epochs. Linear fit of the form $\langle |s^\mu_{j,\theta}| \rangle_\mu = aj + b$ is applied to sine angles across inner iterations $j = 1, \ldots, M$, with $M=10$.}
\centering
\footnotesize
\label{tab:res_super_stats_scientific}
\begin{tabular}{rcccccc}
\toprule
Epoch & $\frac{1}{M}\sum_{j=1}^M \langle|s^\mu_{j,\theta}|\rangle_\mu$ & $\left( \frac{1}{M}\sum_{j=1}^M  \langle |s^\mu_{j,\theta}| \rangle_\mu \right)^M$ & $\Delta_+$ & $\Delta_-$ & $a$ & $b$ \\
\midrule
0   & 0.592766 & $5.36 \times 10^{-3}$ & 0.246203 & -0.518429 & 0.009823 & 0.489628 \\
10  & 0.399244 & $1.03 \times 10^{-4}$  & 0.346335 & -0.119267 & 0.006651 & 0.362663 \\
20  & 0.373864 & $5.30 \times 10^{-5}$   & 0.360928 & -0.148762 & 0.004453 & 0.349371 \\
40  & 0.342721 & $2.20 \times 10^{-5}$   & 0.402198 & -0.127696 & 0.005862 & 0.310482 \\
60  & 0.335581 & $1.80 \times 10^{-5}$   & 0.352728 & -0.122402 & 0.004203 & 0.312462 \\
80  & 0.328283 & $1.50 \times 10^{-5}$   & 0.382090 & -0.114524 & 0.006195 & 0.294211 \\
150 & 0.317659 & $1.00 \times 10^{-5}$   & 0.414758 & -0.111201 & 0.003036 & 0.300960 \\
\bottomrule
\end{tabular}
\end{table}
\begin{figure}[h!]
\centering
\includegraphics[width=0.99\linewidth]{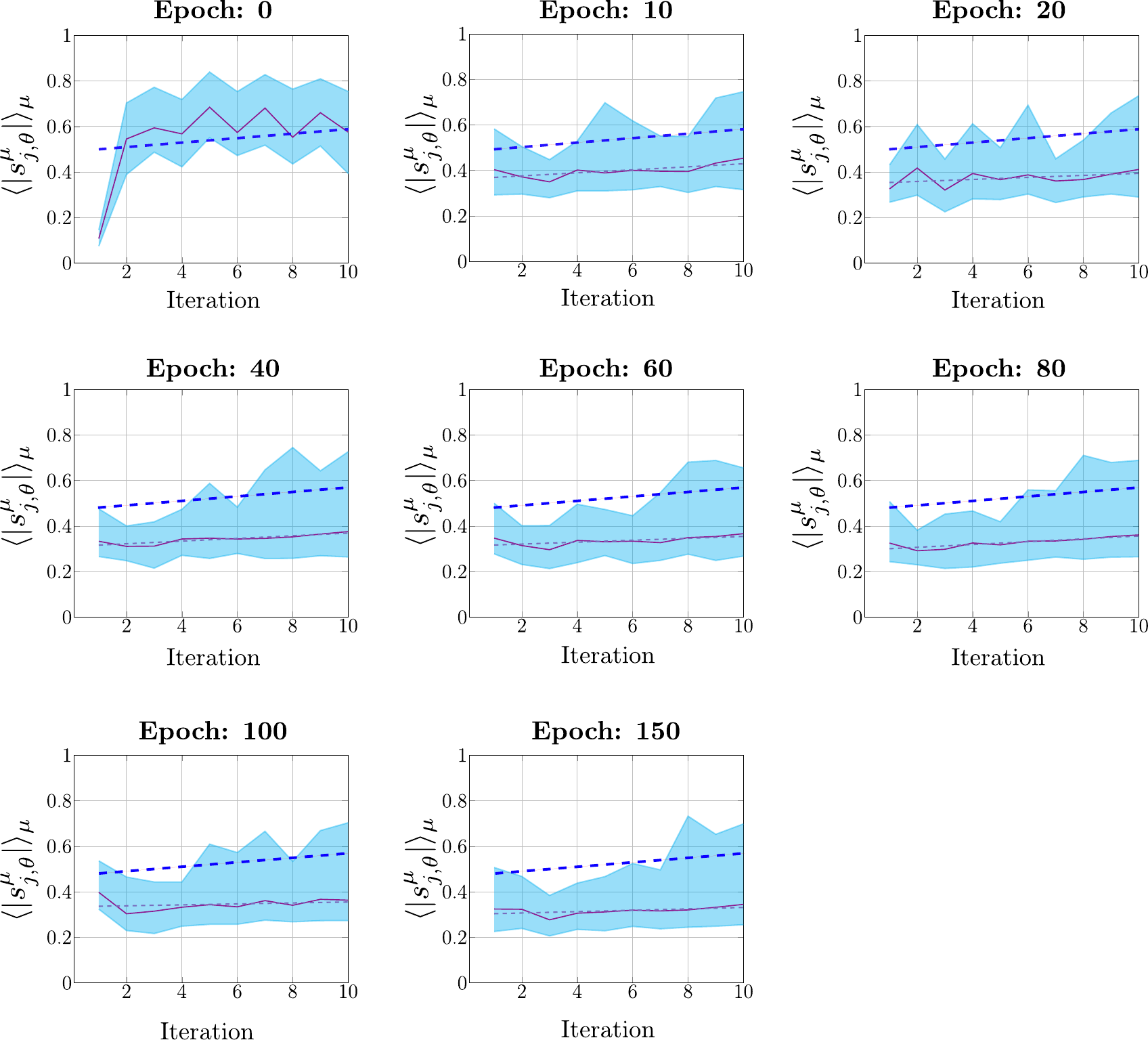} 
\caption{Averaged sine angles $\langle |s^\mu_{j,\theta}| \rangle_\mu$ for $j = 1, \ldots, M$, with $\mu \in \mathcal{P}_\text{validation}$ (solid purple line). The shaded region denotes the envelope of $|s^\mu_{j,\theta}|$ curves across the parameter space. The thick dashed blue line indicates the linear fit of $\langle |s^\mu_{j,\theta}| \rangle_\mu$ at epoch 0, while purple dotted lines show the linear fits at subsequent epochs. Lower sine values reflect improved alignment of residuals with Krylov subspace directions.}
\label{fig:sinus_plot}
\end{figure}
Empirically, during dynamic loss minimization, the averaged principal angles 
$\langle |s^\mu_{j,\theta}| \rangle_\mu$  over the first  $M$ iterations follow an approximately linear trend in the iteration index $j$,
\begin{equation*}
    \left\langle |s^\mu_{j,\theta}| \right\rangle_\mu  \approx a j + b, \quad \text{for } j = 1, \dots, M\, .
\end{equation*}
While the slope $a$ remains nearly constant across training epochs (on the order of $10^{-3}$), the intercept $b$ decreases progressively during optimization, from an initial value of approximately $b \approx 0.48$ down to $b \approx 0.30$ after fine-tuning. Results are summarized in Table~\ref{tab:residual} and Figure~\ref{fig:sinus_plot}.

This decrease reflects improved initial alignment of the residual with the Krylov subspace within the first $M$ iterations. Assuming that $a$ is negligible, we may estimate the product of sines governing the convergence factor via
\begin{equation}
    \prod_{j=1}^M \left\langle |s^\mu_{j,\theta}| \right\rangle_\mu \approx \left( \frac{1}{M} \sum_{j=1}^M \left\langle |s^\mu_{j,\theta}| \right\rangle_\mu \right)^M,
\label{eq:mean_estimate}
\end{equation}
which approximates the asymptotic bound for the Krylov Subspace method convergence. By looking at the estimate above, we can observe a substantial drop (two orders of magnitude) of its value, from $5.3 \times 10^{-3}$  to $1.0 \times 10^{-5}$, experienced during the dynamic training stage. Given the tight relation between principal angles $s_j$ and the convergence rate of the FGMRES algorithm, we expect the same amount of reduction in the value of the relative residual.

\subsubsection{Residual Decay Characterization}

The progressive alignment between the principal angles and the Krylov subspace basis, discussed in Section~\ref{sec:angles}, has a direct impact on the convergence profile of the iterative solver. Indeed, the reduction of  $|s_j|$ induces a contraction of the residual norm, consistent with the geometric recurrence relation detailed in~\eqref{eq:sinus_recurrence}. 

The evolution of the mean residuals can be tracked by a least-squares fitting with an exponential model 
\begin{equation*}
    \left\langle\frac{\| \mathsf{ r^\mu_j} \|}{\| \mathsf{r^\mu_0} \|} \right\rangle_\mu= e^\beta e^{\alpha j}, \quad j = 1,\ldots,M,
\end{equation*}
where $\alpha$ encodes the decay rate and $\beta$ controls the residual offset. The fitting is performed over $M = 10$ iterations, in line with the optimization horizon used in the dynamic loss functional introduced in~\eqref{eq:dynamic_loss}.
At the onset of training (epoch $0$, mean sine $\approx 0.59$), when the network $\mathcal{N}_\theta$ has only undergone static pretraining via the residual-based loss $\mathcal{L}_\text{static}$ (cf. Section~\ref{sec:static_training}), we observe a relative slow decay with parameters $\beta \approx -1.85$ and $\alpha \approx -0.4$  which leads to an average relative residual of $\sim 10^{-3}$ after 10 iterations. As fine tuning progresses, the minimization of the dynamic loss $\mathcal{L}_{\text{dynamic}(M)}$ actively promotes alignment of the residual with the Krylov subspace; thus $\alpha$ decreases significantly, from $\alpha \approx -0.4$ at epoch $0$ to $\alpha \approx -1.13$ at epoch $150$. The reduction in $\alpha$ reflects an acceleration in the asymptotic convergence rate.
\begin{table}[h!]
\label{tab:residual}
\caption{Evolution of relative residual at iteration M. Exponential fit of the form $\langle ||\mathsf{r_{j}^\mu||/||r^\mu_0}||\rangle_\mu = e^{\beta}e^{\alpha j}$ is applied to relative residuals for inner iterations $j = 1, \ldots, M$.} \vspace{10pt}
\centering
\footnotesize
\label{tab:res_super_stats_sci}
\begin{tabular}{r@{\hskip 30pt} c@{\hskip 30pt}c@{\hskip 30pt}c@{\hskip 30pt}c@{\hskip 30pt}c}
\toprule
Epoch &  $\langle ||\mathsf{r_{M}^\mu||/||r^\mu_0}||\rangle_\mu$& $\Delta_+$& $\Delta_-$&  $\alpha$ & $\beta$ \\
\midrule
0   & $1.20 \times 10^{-3}$ & 0.003509 & -0.000849 & -0.489807 & -1.851620 \\
10  & $1.05 \times 10^{-4}$ & 0.000415 & -0.000079 & -0.919889 & -0.123565 \\
20  & $5.50 \times 10^{-5}$ & 0.000332 & -0.000041 & -0.975239 & -0.161445 \\
40  & $2.70 \times 10^{-5}$ & 0.000257 & -0.000021 & -1.049013 & -0.225478 \\
60  & $2.30 \times 10^{-5}$ & 0.000252 & -0.000018 & -1.077473 & -0.136004 \\
80  & $1.80 \times 10^{-5}$ & 0.000298 & -0.000015 & -1.090334 & -0.240775 \\
100 & $2.90 \times 10^{-5}$ & 0.000274 & -0.000024 & -1.057390 & -0.052897 \\
150 & $1.30 \times 10^{-5}$ & 0.000209 & -0.000010 & -1.134431 & -0.095151 \\
\bottomrule
\end{tabular}
\end{table}
\begin{figure}[h!]
\centering
\includegraphics[width=0.99\linewidth]{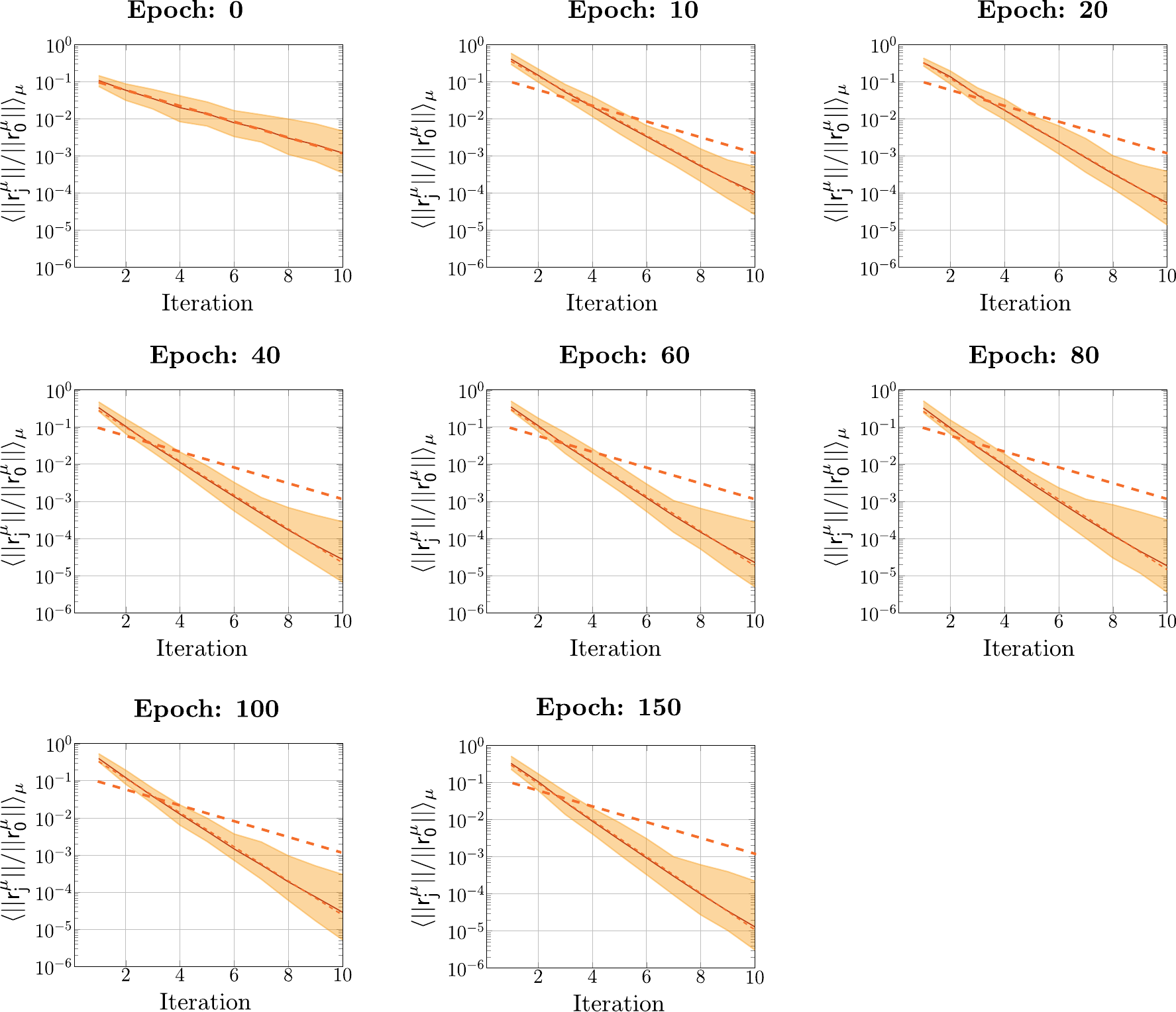}
\caption{Averaged relative residuals $\langle ||\mathsf{r_{j}^\mu||/||r^\mu_0}||\rangle_\mu$ for $j = 1, \ldots, M$, with $\mu \in \mathcal{P}_\text{validation}$ (solid red line). The shaded region denotes the envelope of $||\mathsf{r_{j}^\mu||/||r^\mu_0}||$ curves across the parameter space. The thick dashed red line indicates the linear fit of $\langle ||\mathsf{r_{j}^\mu||/||r^\mu_0}||\rangle_\mu$  at epoch 0, while red dotted lines show the exponential fits at subsequent epochs.}
\label{fig:res_plot}
\end{figure}

For $M=10$, the mean relative residual (averaged over the validation batch) decreases from $1.12 \times 10^{-3}$ at epoch $0$ to $1.3 \times 10^{-5}$ at the final epoch. This represents a two-order-of-magnitude improvement in the reduction of the residual norm achieved through the fine-tuning phase. 

As expected, these results perfectly align with the estimates in Eq.~\ref{eq:mean_estimate} reported in Table \ref{tab:sinus}.
All results are summarized in Table~\ref{tab:residual} and Figure~\ref{fig:res_plot}, which report the fitted parameters and associated statistics.

\subsection{Localized Control over Krylov Dynamics}

An important observation is that the dynamic loss structure defined by $\mathcal{L}_{\text{dynamic}(M)}$ induces a localized improvement in Krylov geometry: by restricting the supervision window to $M = 10$ Arnoldi iterations, the neural preconditioner $\mathcal{N}^\text{dy}_\theta$ \textit{specializes} in optimizing early-stage convergence—where computational savings are most impactful. Beyond the supervised horizon ($j > M$), the principal angles $|s^\mu_{j,\theta}|$ tend to increase, reflecting the absence of penalization in later iterations.

Crucially, this behavior is not a limitation but a controllable feature: by tuning the value of $M$, one can trade depth for early accuracy, thereby aligning the optimization objective with specific solver performance targets. This targeted control provides a flexible and interpretable framework for tailoring the action of learned preconditioners to varying computational demands.

\begin{figure}[h]
\centering
\includegraphics[width=0.99\linewidth]{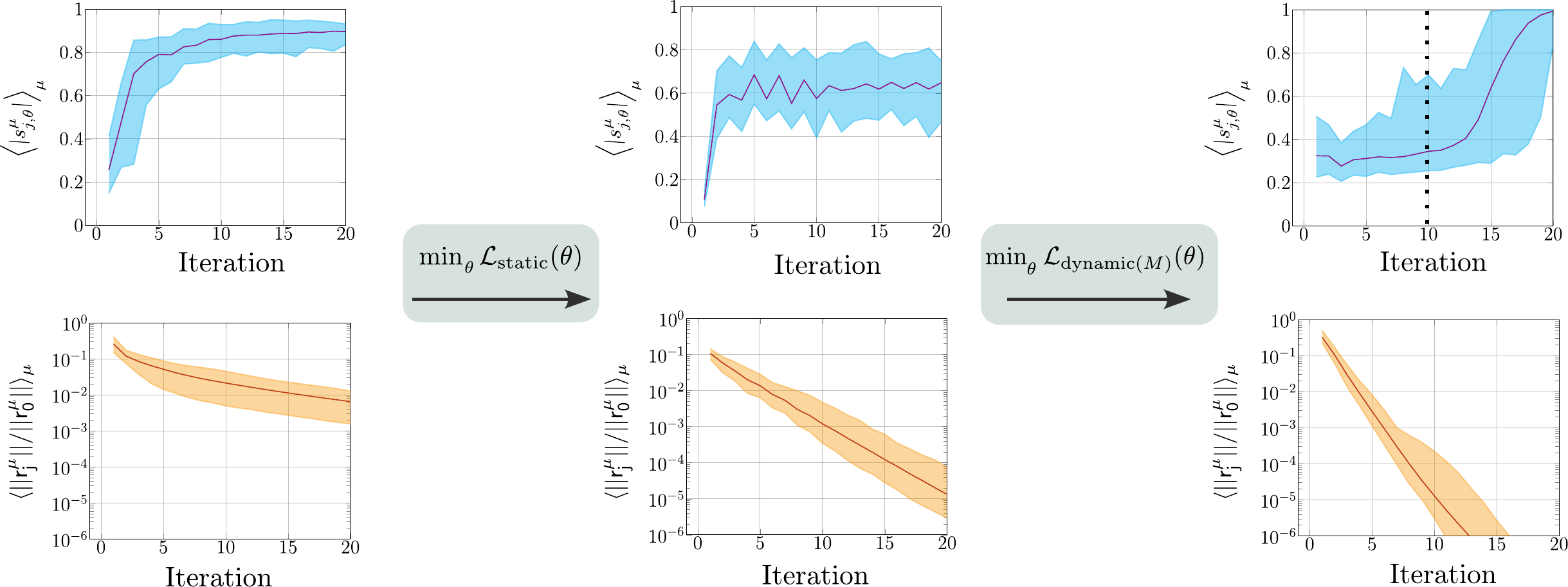}
\caption{Comparison of mean principal angles $\langle |s^\mu_{j,\theta}| \rangle_\mu$ and relative residuals $\langle\mathsf{||r_{j}^\mu||/||r^\mu_0}||\rangle_\mu$ for $j = 1, \ldots, 2M$, after static training via $\mathcal{L}_\text{static}(\theta)$ (\textbf{left}) and after dynamic fine-tuning via $\mathcal{L}_{\text{dynamic}(M)}(\theta)$ with $M=10$ (\textbf{right}).}
\label{fig:sin_res_combined}
\end{figure}

\subsection{Iteration Count Analysis}

We now evaluate the performance of the dynamically fine-tuned neural preconditioner $\mathcal{N}_\theta^\text{dy}$ obtained via the two-stage training protocol outlined in Section~\ref{sec:training}. The benchmark involves a parametric family of coupled 3D–1D problems, where the 1D domain represents vascular networks modeled as randomly generated graphs of different geometric complexity. This setting introduces significant variability in the system matrices, making it a robust test case for generalization across operator structures.

The objective of the preconditioner is to improve both spectral and geometric properties of the system matrix to accelerate Krylov subspace convergence. To assess this, we measure the average number of FGMRES iterations required to reduce the relative residual norm below $10^{-6}$ across a test set of $100$ previously unseen graph configurations. Without preconditioning, the solver requires, on average $147.5$ iterations, with values ranging from $100$ to $245$. With the statically trained preconditioner, the count drops sharply to $26.71$, with a narrower range of $23$ to $31$ iterations. Following dynamic fine-tuning, the average count further decreases to $12.86$, with values ranging from $11$ to $16$.

These results show that the neural preconditioner $\mathcal{N}_\theta$ reduces the average iteration count by more than an order of magnitude across the parametric problem family $(\mathsf{A^\mu}, \mathsf{b^\mu})$. The additional gain from dynamic tuning confirms that the Krylov subspace generation process itself is being effectively manipulated by the neural preconditioner.

\begin{figure}
    \centering
    \includegraphics[width=1.0\linewidth]{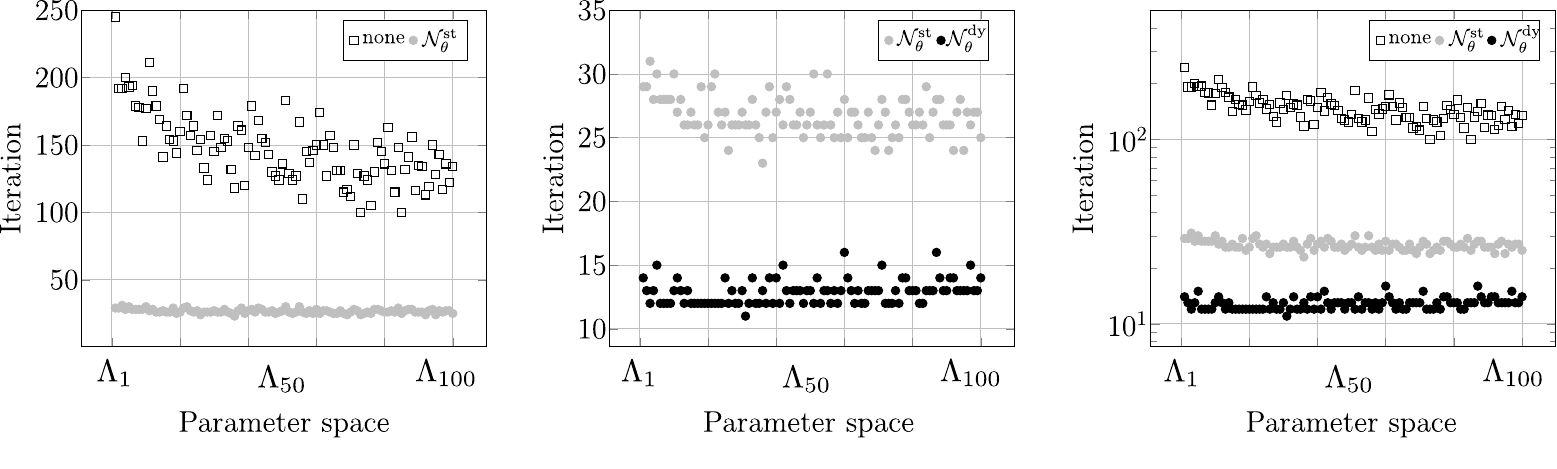}
    \caption{FGMRES iteration counts across 100 test 1D graph instances $\Lambda_i$. \textbf{Left:} unpreconditioned vs.\ statically trained $\mathcal{N}_\theta^\text{st}$. \textbf{Center:} static vs.\ dynamically fine-tuned $\mathcal{N}_\theta^\text{dy}$. \textbf{Right:} global overview across all configurations.}
\end{figure}

\section{Conclusion and Perspectives}

This work introduces a novel training strategy for neural preconditioners that leverages the geometric structure of Krylov subspace methods, specifically the Generalized Minimal Residual (GMRES) algorithm. The core contribution lies in the formulation of a dynamic fine-tuning phase, wherein the loss functional directly optimizes the subspace angles \( s_{i,\theta} \) that govern the convergence behavior of iterative solvers. By embedding this geometric insight into the training process, we offer a principled and solver-aligned approach to performance enhancement.

A central advantage of our method is that it preserves the strengths of unsupervised learning—such as straightforward data generation and independence from ground truth solutions—while simultaneously introducing a transparent and solver-integrated performance metric. The use of a differentiable formulation of the Flexible GMRES algorithm enables efficient gradient-based optimization through backpropagation across the solver's iterative process. The dynamic residual-based loss we propose captures the alignment between the residual vector and the Krylov subspace at each iteration, quantified by the sine of their principal angle, thereby offering a geometrically meaningful proxy for convergence rate.

The effectiveness of this solver-aware training scheme is validated by substantial reductions in iteration counts for parameter-dependent linear systems derived from mixed-dimensional partial differential equations (PDEs)—a class of problems often characterized by ill-conditioning and geometric complexity. Building upon our earlier work on unsupervised neural preconditioners for mixed-dimensional models, we demonstrate that the combined static-dynamic training strategy significantly enhances the preconditioner’s adaptability to heterogeneous problem instances. In our numerical tests, the average number of FGMRES iterations dropped from 147.5 (unpreconditioned case) to 26.71 after static training, and further to 12.86 following dynamic fine-tuning. This outcome attests to the ability of the neural preconditioner to influence the subspace geometry and accelerate solver convergence directly.

Despite these encouraging results, several important directions remain for future research. Extending the proposed approach to alternative iterative solvers and different classes of preconditioners represents a natural next step. Investigating hybrid training strategies that blend static and dynamic loss formulations, or integrate additional physics-informed components, could yield further improvements in robustness and generalization. A particularly relevant challenge is to assess the preconditioner's performance on out-of-distribution samples, including unseen geometries and parametric regimes, to better understand its extrapolative capabilities.

From a theoretical standpoint, formalizing the role of spectral data augmentation and characterizing the learning dynamics of Krylov-aware neural networks are crucial for building a more rigorous foundation for learning-based preconditioning. Moreover, current limitations related to the use of convolutional neural networks on structured grids suggest the need for architectural generalizations. Mesh-informed neural networks \cite{FrancoMINN} or graph-based models  \cite{gnns} could provide the necessary flexibility to extend applicability to unstructured meshes and more general computational domains.

Finally, scaling up the approach to large-scale problems remains a critical avenue. Optimizing multi-resolution network design, exploiting GPU-parallel computations, and integrating ensemble methods can substantially improve throughput in multi-query settings such as parameter studies, control, or uncertainty quantification.

In conclusion, this study establishes the feasibility and effectiveness of neural preconditioning guided by Krylov subspace geometry as a powerful tool for accelerating iterative solvers in complex PDE applications. It provides a promising foundation for advancing the integration of deep learning and numerical linear algebra in scientific computing.

\section*{Acknowledgments}
PZ acknowledges the support of the MUR PRIN 2022 grant No. 2022WKWZA8 \emph{Immersed methods for multiscale and multiphysics problems} (IMMEDIATE) part of the Next Generation EU program Mission 4, Comp. 2, CUP D53D23006010006. The present research is part of the activities of the \emph{Dipartimento di Eccellenza} 2023-2027, Department of Mathematics, Politecnico di Milano.
AC acknowledges the support of the MUR PRIN 2022 project \emph{Evolution problems involving interacting scales} project code 2022M9BKBC, Grant No. CUP  D53D23005880006.
All authors are members of the Gruppo Nazionale per il Calcolo Scientifico (GNCS), Istituto Nazionale di Alta Matematica (INdAM).

\bibliography{references}

\end{document}